\documentclass{amsart}
\usepackage{amsmath,cite}
\usepackage{mathtools}
\usepackage{color}
\usepackage[colorlinks=true,urlcolor=blue,
citecolor=red,linkcolor=blue,linktocpage,pdfpagelabels,bookmarksnumbered,bookmarksopen]{hyperref}
\usepackage[english]{babel}

\numberwithin{equation}{section}
\newtheorem{theorem}{Theorem}[section]
\newtheorem{proposition}[theorem]{Proposition}

\theoremstyle{definition}

\newcommand{\R}{{\mathbb R}}
\renewcommand{\O}{{\mathcal O}}
\newcommand{\iu}{{\rm i}}
\newcommand{\eps}{\varepsilon}
\newcommand{\C}{{\mathbb C}}

\DeclarePairedDelimiter{\nhu}{\lVert}{\rVert_{H^1(\R^N)}}

\DeclarePairedDelimiter{\nldd}{\lVert}{\rVert^2_{L^2(\R^N)}}
\DeclarePairedDelimiter{\nlpsps}{\lVert}{\rVert^{p+1}_{L^{p+1}(\R^N)}}

\title[On the stability of Klein-Gordon 
equations]{On the stability of standing waves of Klein-Gordon equations in a semiclassical regime}

\author[M.~Ghimenti]{Marco Ghimenti}
\author[S.~Le Coz]{Stefan Le Coz}
\author[M.~Squassina]{Marco Squassina}

\address{Dipartimento di Matematica Applicata
\newline\indent
Universit\`a di Pisa
\newline\indent
Via F. Buonarroti 1/c, 56127 Pisa
\newline\indent
Italia}
\email{ghimenti@mail.dm.unipi.it}

\address{Institut de Math\'ematiques de Toulouse,
\newline\indent
Universit\'e Paul Sabatier
\newline\indent
118 route de Narbonne, 31062 Toulouse Cedex 9
\newline\indent
France}
\email{slecoz@math.univ-toulouse.fr}

\address{Dipartimento di Informatica 
\newline\indent
Universit\`a di Verona
\newline\indent
Strada Le Grazie 15, 37134 Verona
\newline\indent
Italia}
\email{marco.squassina@univr.it}

\thanks{
Research supported in part by 2009 MIUR project: {\em ``Variational and Topological
Methods in the Study of Nonlinear Phenomena''},  ANR project ESONSE and 
GNAMPA project {\em ``Dinamica e propriet\`a di soluzioni concentrate in teorie di campo nonlineari'' }
}

\date{\today}

\begin{document}

\begin{abstract}
We investigate the orbital stability and instability of standing waves for two classes of Klein-Gordon
equations in the semi-classical regime.
\end{abstract}

\maketitle

\section{Introduction and results}

The nonlinear Klein-Gordon equation 
\begin{equation}
	\label{kg}
\varepsilon^2 u_{tt}-\eps^2\Delta u+mu-|u|^{p-1}u=0 \qquad (t,x)\in\R^+\times \R^N,
\end{equation}
where $\eps,m>0$, $p>1$ for $N=1,2$ and $1<p<(N+2)/(N-2)$ for $N\geq 3$, is one of the simplest nonlinear 
partial differential equations invariant for the Poincar\'e group.
We are interested in the study of the nonlinear Klein Gordon equation in presence of a potential  
depending on the space variable.
Two different choices are viable. We can simply add a potential term $W(x)u$ to equation \eqref{kg}.
This case has been studied, for the linear wave equation, for example, by Beals and Strauss in \cite{BS}. 
Otherwise, if we want to fully preserve the invariance for the Poincar\'e group of \eqref{kg}, 
we have to change the temporal derivative $\eps^2\partial_{tt}$ with a covariant derivative, depending on the potential 
$D^2_{tt}$, where $D_t= \eps\partial_t +\iu V(x)$. This second approach is classical, in quantum electrodynamics, when considering 
electromagnetic waves. The first approach leads us to consider the equation 
\begin{equation}
\label{seconmod}
\eps^2 u_{tt}-\eps^2\Delta u+mu-Wu-|u|^{p-1}u=0, \qquad\text{in }\R^N, 
\end{equation}
while the second one to investigate the equation 
\begin{equation}
\label{primod}
\eps^2 u_{tt}+2\iu\eps Vu_t-\eps^2\Delta u+mu-V^2u-|u|^{p-1}u=0, \qquad\text{in }\R^N.
\end{equation}
In this paper, we treat simultaneously the two previous Klein-Gordon equations by studying 
\begin{equation}
	\label{genequation}
\eps^2 u_{tt}+2\iu\eps Vu_t-\eps^2\Delta u+mu-Wu-|u|^{p-1}u=0, \qquad \text{ in }\R^N,
\end{equation}
where $u:\R^N\times\R\to\C$, $\eps>0$ and $V,W$ are real valued potential functions.
Equation \eqref{genequation} formally yields to \eqref{primod} for the choice $W=V^2$ as well as to \eqref{seconmod} when $V=0$.
We shall study the stability of standing waves of this equation in the semiclassical regime $\eps\to 0$.
It admits standing waves, namely solutions of the form $u(x,t)=e^{\iu\omega t/\eps}\varphi_\omega(x/\eps)$, 
where $\omega\in\R$ and $\varphi_\omega$ satisfies
\begin{equation}
	\label{ellip-omega}
-\Delta \varphi_\omega+\big(m-\omega^2-2\omega V(\eps y)-W(\eps y)\big)\varphi_\omega-|\varphi_\omega|^{p-1}\varphi_\omega=0, \qquad\text{in }\R^N.
\end{equation}
To ensure existence of solutions to \eqref{ellip-omega} for $\eps$ close to $0$, we assume the following. Let $V$ and $W$  be $\mathcal{C}^2$. For the function
$$
Z(y)=m-\omega^2-2\omega V(y)-W(y),\,\,\ y\in\R^N
$$
there exists $x_0\in\R^N$
such that 
\begin{equation}
	\label{criticality}
\nabla Z(x_0) = 0,\qquad \nabla^2 Z(x_0)\text{ is non-degenerate.}
\end{equation}
Furthermore, we assume that
\begin{equation}
	\label{positivity}
Z(x_0)=m-\omega^2-2\omega V(x_0)-W(x_0)>0.
\end{equation}

Under these hypotheses, it is well-know (see e.g. \cite{AmBaCi97} or \cite[Section 8.2]{AmMa06}) that when $\eps$ is close to $0$ the equation
\eqref{ellip-omega} admits a family of positive, exponentially decaying, solutions $\varphi_\omega\subset H^1(\R^N)$ (hiding the dependance 
upon $\eps$). More precisely, there exists $\xi_\eps\in\R^N$ and $\psi_\omega\in H^1(\R^N)$ such that $\varphi_{\omega}(\cdot)= \psi_\omega(\cdot-\xi_\eps)+\O(\eps^2)$ in $H^1(\R^N)$
as $\eps\to 0$, where $\xi_\eps=o(\eps)$ and
\begin{equation*}
-\Delta \psi_\omega+Z(x_0)\psi_\omega=|\psi_\omega|^{p-1}\psi_\omega, \qquad\text{in $\R^N$}.
\end{equation*}

We are interested in the (orbital) stability or instability of standing waves when $\eps$ goes to $0$.  

A standing wave of \eqref{genequation} is said to be \emph{(orbitally) stable} if any solution of \eqref{genequation} 
starting close to the standing wave remains close for all time, up to the invariances 
of the equation. More precisely, we say that $e^{\frac{\iu\omega t}{\eps}}\varphi_\omega\left(\frac{x}{\eps}\right)$ 
is stable if for all $\eta>0$ there exists $\delta>0$ such that for all $u_0\in H^1(\R^N)$ verifying $\nhu{u_0-\varphi_\omega}<\delta$ 
the solution $u(x,t)$ of \eqref{genequation} with initial data $u_0$ satisfies 
\[
\sup_{t\in\R}\inf_{\theta\in\R}\nhu{u-e^{i\theta}\varphi_\omega}<\eta.
\]
Since the pioneering works \cite{BeCa81,CaLi82,GrShSt87,GrShSt90,We83,We85}, the study of orbital stability for standing waves of 
dispersive PDE has attracted a lot of attention. Among many others, one can refer to \cite{JeLe06,JeLe09,LeFuFiKsSi08}; see also 
the books and surveys \cite{Ca03,Le09,St08,Ta09} and the references therein. Relatively few works \cite{IaLe09,LiWe08,Oh89} are concerned 
with stability at the semi-classical limit for Schr\"odinger type equations. For Klein-Gordon equations, after the ground works \cite{ShSt85,Sh85}, 
there has been a recent interest for stability by blow-up \cite{OhTo05,OhTo06,OhTo07,LiOhTo07}. 

To study stability, we first rewrite \eqref{genequation} in Hamiltonian form
\begin{equation}\label{eq:hamilt}
\eps\frac{\partial U}{\partial t}=JE'(U),
\end{equation}
where $U=\begin{pmatrix}u\\v\end{pmatrix}$, $J=\begin{pmatrix}0 & 1\\-1 & 0\end{pmatrix}$, and 
\begin{multline*}
E(U)=\frac12\nldd{v-i Vu}+\frac12\nldd{\nabla u}+m\frac12\nldd{u}\\-\frac12\int_{\R^N}W|u|^2dx-\frac{1}{p+1}\nlpsps{u}.
\end{multline*}
It is easy to see that if $u$ solves \eqref{genequation} and $v$ is defined by $v:=\eps u_t+i V u$, then $U=\begin{pmatrix}u\\v\end{pmatrix}$ solves \eqref{eq:hamilt}.
The charge $Q(U)$ is defined by
\[
Q(U)=\Im\int_{\R^N}\bar{u}vdx.
\]
In particular, for a standing wave $u=e^{i\omega t/\eps}\varphi_\omega(x/\eps)$, the charge is given by
\begin{equation}\label{eq:defQ}
Q(\varphi_\omega):=Q(U)=\eps^N\left(\omega\nldd{\varphi_\omega}+\int_{\R^N} V(\eps y)|\varphi_\omega|^2\right).
\end{equation}

According to the theory developed in \cite{GrShSt87,GrShSt90}, a standing wave $e^{\frac{\iu\omega t}{\eps}}\varphi_\omega\left(\frac{x}{\eps}\right)$ is stable if two conditions are satisfied. 
\begin{itemize}
\item[(i)] The \emph{Slope Condition}: $\displaystyle \frac{\partial}{\partial \omega}Q(\varphi_\omega)<0$.
\item[(ii)] The \emph{Spectral Condition}: $L_\eps:=-\Delta+Z(\eps y)-p|\varphi_\omega|^{p-1}$ has exactly one negative eigenvalue and is non degenerate.
\end{itemize}
On the other hand, denote by $n(L_\eps)$ the number of negative eigenvalues of $L_\eps$ and set $p(\omega) =
0$ if $\displaystyle \frac{\partial}{\partial \omega}Q(\varphi_\omega)>0$, $p(\omega) = 1$ if $\displaystyle \frac{\partial}{\partial \omega}Q(\varphi_\omega)<0$. 
Then the standing wave is unstable if $n(L_\eps)-p(\omega)$ is odd.

Then we have the following
\begin{theorem}
	\label{mainthm}
Assume that conditions \eqref{criticality}-\eqref{positivity} hold. Then, we have the following facts.
\begin{enumerate}
\item 
If $p<1+4/N$, then the Slope Condition $\frac{\partial}{\partial\omega} Q(\varphi_\omega)<0$ is fulfilled if
\[
Z(x_0)<(\omega+V(x_0))^2\Big(\frac{4}{p-1}-N\Big) \qquad\text{(non-critical case)}
\]
or if 
\[
\left\{
\begin{array}{l}
Z(x_0)=(\omega+V(x_0))^2\Big(\frac{4}{p-1}-N\Big), \\
\left(\Delta Z(0)-\Delta V(0)\left(1+\frac{2(\omega+V(0))}{Z(0)}\right)\right)<0,
\end{array}
\right.
\qquad\text{(critical case).}
\]
\item
If $p\geq 1+4/N$, then we always have $\frac{\partial}{\partial\omega} Q(\varphi_\omega)>0$.
\item
We have the equality $n(L_\eps)=n(\nabla^2Z(x_0))+1$, where $n(\nabla^2Z(x_0))$ is the number of negative eigenvalues of $\nabla^2Z(x_0)$.
\end{enumerate}
In particular, the standing waves $e^{i\omega t}\varphi_\omega$ are stable if $x_0$ is non-degenerate local minimum of $Z$, $p<1+4/N$ and 
\[
Z(x_0)<(\omega+V(x_0))^2\Big(\frac{4}{p-1}-N\Big).
\]
\end{theorem}

Note that, conversely to what was happening in the case of Schr\"odinger equations studied in \cite{LiWe08}, the values of the potentials $V$ and $W$ 
at $x_0$ comes into play for the Slope Condition even in the noncritical case. Note also that only the local behavior of $Z$ around $x_0$ influences the stability or instability.

Notations : Most of the objects we consider will depend both on $\eps$ and $\omega$. We will emphasize the most 
important parameter by indicating it as a subscript, the dependence in the other parameter being understood.

\section{Proof of Theorem~\ref{mainthm}}

In this section, we prove Theorem \ref{mainthm}. We start be focusing on the 
Slope Condition and then we study the Spectral Condition. For the sake of simplicity in notations and without loss of generality, in the rest of this section we assume that $x_0=0$.

\subsection{The Slope Condition}

We start with the noncritical case.

\subsubsection{Noncritical case}

We assume that 
\[
Z(0)\neq(\omega+V(0))^2\Big(\frac{4}{p-1}-N\Big).
\]

We first rewrite $Q(\varphi_\omega)$ by expanding $V(\eps y)$ and using the exponential decay of $\varphi_\omega$:
\[
Q(\varphi_\omega)=\eps^N\left( \omega+V(0) \right)\nldd{\varphi_\omega}+\O(\eps^{N+1}).
\]
Therefore, since
\begin{equation}
	\label{Q1}
\frac{\partial}{\partial\omega} Q(\varphi_\omega)=\eps^N
\|\varphi_\omega\|_{L^2(\R^N)}^2+\eps^N(\omega+V(0) )\frac{\partial}{\partial\omega}\|\varphi_\omega\|_{L^2(\R^N)}^2+\O(\eps^{N+1}),
\end{equation}
to evaluate the sign of the map
$\omega\mapsto \frac{\partial}{\partial\omega} Q(\varphi_\omega)$
one should compute the quantity
\begin{equation}
	\label{Q2}
\frac{\partial}{\partial\omega}\|\varphi_\omega\|_{L^2(\R^N)}^2
=2\int_{\R^N} R_\omega\varphi_\omega,
\end{equation}
where $R_\omega(x):=\frac{\partial\varphi_\omega}{\partial\omega}(x)$.

We remark that differentiation of \eqref{ellip-omega} with respect to $\omega$ easily yields 
\begin{equation}
	\label{identitt}
L_\eps R_\omega=2(\omega +V(\eps y))\varphi_\omega.
\end{equation}
If we now introduce the rescaling $\varphi_\omega(x)=\lambda^{\frac{1}{p-1}}\varphi_\lambda(\sqrt{\lambda}x)$,
it follows that $\varphi_\lambda$ satisfies 
\begin{equation}
	\label{ellip-lambda}
-\Delta \varphi_\lambda+
\lambda^{-1}Z\left(\frac{\eps y}{\sqrt{\lambda}}\right)\varphi_\lambda
-|\varphi_\lambda|^{p-1}\varphi_\lambda=0, \qquad\text{in $\R^N$}.
\end{equation}
Now, differentiating equation \eqref{ellip-lambda} with respect to $\lambda$ and denoting
$T_\lambda=\frac{\partial\varphi_\lambda}{\partial\lambda}_{|\lambda=1}$ yields 
\begin{equation}
	\label{oplin}
L_\eps T_\lambda-
Z(\eps y)\varphi_\omega
-\frac{1}{2} \eps y\cdot\nabla Z(\eps y)\varphi_\omega=0.
\end{equation}

Since $0$ is a critical point of $Z$, a Taylor expansion gives
\begin{gather}
	\label{sviluppo1}
Z(\eps y)=Z(0)+\O(\eps^2|y|^2),\\
	\label{sviluppo2}
\frac{1}{2} \eps y\cdot\nabla Z(\eps y)=\O(\eps^2|y|^2).
\end{gather}
Then, from \eqref{oplin},  as $\eps\to 0$ we have
\begin{equation}
	\label{oplin-appr}
L_\eps T_\lambda=Z(0)\varphi_\omega+\O(\eps^2|y^2|\varphi_\omega),\quad\text{in $\R^N$}.
\end{equation}
Then, in turn, taking into account identity~\eqref{identitt}  we get
\begin{align}
Z(0)\int_{\R^N} R_\omega\varphi_\omega&=\int_{\R^N} R_\omega L_\eps T_\lambda+\O(\eps^2) \nonumber\\
&=\int_{\R^N} L_\eps R_\omega T_\lambda+\O(\eps^2) \nonumber \\
&=\int_{\R^N} 2(\omega +V(\eps y))\varphi_\omega T_\lambda+\O(\eps^2) \label{eq:nondegcase}\\
&=2(\omega+V(0))\int_{\R^N}\varphi_\omega T_\lambda+\O(\eps)  \nonumber\\
&=(\omega+V(0))\frac{\partial}{\partial\lambda}{\|\varphi_\lambda\|_{L^2(\R^N)}^2}_{|\lambda=1}+\O(\eps)  \nonumber\\
&=(\omega+V(0))\Big(\frac{N}{2}-\frac{2}{p-1}\Big)\|\varphi_\omega\|_{L^2(\R^N)}^2+\O(\eps). \nonumber
\end{align}
In conclusion, by combining \eqref{Q1}, \eqref{Q2} and \eqref{oplin-appr}, we have
\[
\frac{\partial}{\partial\omega} Q(\varphi_\omega)=\eps^N
\left(1+\frac{(\omega+V(0))^2}{Z(0)}\Big(N-\frac{4}{p-1}\Big)\right)
\|\varphi_\omega\|_{L^2(\R^N)}^2+\O(\eps^{N+1}).
\]
Then, taking into account the fact that $Z(0)>0$ and that $\varphi_\omega$ converges to $\psi_\omega$ in $L^2(\R^N)$ as $\eps\to 0$, 
the sign of $\frac{\partial}{\partial\omega} Q(\varphi_\omega)$ is the sign of 
\[
Z(0)-(\omega+V(0))^2\Big(\frac{4}{p-1}-N\Big).
\]

\subsubsection{Critical case}
We assume now that 
\begin{equation}\label{eq:degencase}
Z(0)=(\omega+V(0))^2\Big(\frac{4}{p-1}-N\Big).
\end{equation}
In the critical case, the term of order $\eps^N$ in the 
expansion of $\frac{\partial }{\partial \omega}Q(\varphi_\omega)$ 
vanishes and we need to calculate the expansion at a higher order.
We first refine \eqref{sviluppo1}-\eqref{sviluppo2}.
\begin{gather*}
Z(\eps y)=Z(0)+\frac{\eps^2}{2}\nabla^2Z(0)(y,y)+\O(\eps^3|y|^3)\\
\frac{1}{2} \eps y\cdot\nabla Z(\eps y)=\frac{\eps^2}{2} \nabla^2 Z(0)(y,y)+\O(\eps^3|y|^3).
\end{gather*}
Then \eqref{oplin} gives
\[
L_\eps T_\lambda=Z(0)\varphi_\omega+\eps^2\nabla^2 Z(0)(y,y)\varphi_\omega+\O(\eps^3|y^3|)\varphi_\omega.
\]
Now, we have 
\begin{equation}\label{eq:hessian1}
Z(0)\int_{\R^N}R_\omega\varphi_\omega=\int_{\R^N}R_\omega L_\eps T_\lambda-\eps^2\int_{\R^N}\nabla^2 Z(0)(y,y)R_\omega\varphi_\omega+\O(\eps^3).
\end{equation}
From \eqref{identitt}, we obtain
\begin{equation}\label{eq:comeback}
\int_{\R^N}R_\omega L_\eps T_\lambda=\int_{\R^N}L_\eps R_\omega T_\lambda=\int_{\R^N} 2(\omega+V(\eps y))\varphi_\omega T_\lambda.
\end{equation}
Expanding the potential $V$ we get
\begin{multline}\label{eq:Vexpansion}
\int_{\R^N} 2V(\eps y)\varphi_\omega T_\lambda=\int_{\R^N} 2V(0)\varphi_\omega T_\lambda+2\eps\int_{\R^N}  
y\cdot\nabla V(0)\varphi_\omega T_\lambda\\+\eps^2\int_{\R^N} \nabla^2V(0)(y,y)\varphi_\omega T_\lambda+\O(\eps^3).
\end{multline}
Note that since $\varphi_\omega=\psi_\omega(\cdot-\xi_\eps)+\O(\eps^2)$ and $\xi_\eps=o(\eps)$,
we have
\begin{equation}\label{eq:cancelation}
2\eps\int_{\R^N}  y\cdot\nabla V(0)\varphi_\omega T_\lambda=2\eps\int_{\R^N}  y\cdot\nabla V(0)\psi_\omega T_\lambda+o(\eps^2)=o(\eps^2)
\end{equation}
where the last cancellation comes from the fact that $\psi_\omega$ is radial.
Coming back to \eqref{eq:comeback} and as in \eqref{eq:nondegcase}, we have 
\begin{multline}\label{eq:hessian2}
\int_{\R^N}R_\omega L_\eps T_\lambda
=(\omega+V(0))\Big(\frac{N}{2}-\frac{2}{p-1}\Big)\|\varphi_\omega\|_{L^2(\R^N)}^2\\+\eps^2\int_{\R^N} \nabla^2V(0)(y,y)\varphi_\omega T_\lambda+o(\eps^2).
\end{multline}

It remains to compute the integrals involving the Hessians in \eqref{eq:hessian1} and \eqref{eq:hessian2}. Since our problem is invariant by an 
orthonormal transformation, we can assume without loss of generality 
that $\nabla^2V(0)=\mathrm{diag}(b_1,\dots,b_N)$. Hence $\nabla^2V(0)(y,y)=\sum_{j=1}^Nb_jy_j^2$. 
Recall also that $T_\lambda$ can be computed explicitly to have 
$$
T_\lambda=-\frac{1}{p-1}\varphi_\omega-\frac{1}{2}y\cdot\nabla\varphi_\omega.
$$ 
Therefore, 
\[
\int_{\R^N} b_jy_j^2\varphi_\omega T_\lambda=-\frac{b_j}{p-1}\int_{\R^N}y_j^2\varphi_\omega^2-\frac{b_j}{2}\sum_{k=1}^N\int_{\R^N}y_j^2y_k\varphi_\omega\frac{\partial\varphi_\omega}{\partial y_k}.
 \]
We have after integration by parts
\[
2\sum_{k=1}^N\int_{\R^N}y_j^2y_k\varphi_\omega\frac{\partial\varphi_\omega}{\partial y_k}=-\sum_{k=1}^N\int_{\R^N}(y_j^2+2\delta_{jk}y_j^2)\varphi_\omega^2=-(N+2)\int_{\R^N}y_j^2\varphi_\omega^2.
\]
Thus
\[
\int_{\R^N} \nabla^2V(0)(y,y)\varphi_\omega T_\lambda=\sum_{j=1}^N\int_{\R^N} b_jy_j^2\varphi_\omega T_\lambda=-\left(\frac{1}{p-1}-\frac{N+2}{4}\right)\sum_{j=1}^Nb_j\int_{\R^N}y_j^2\varphi_\omega^2.
\]
Recall the following expansion in $\eps$ for $R_\omega$ and $\varphi_\omega$.
\[
\varphi_\omega=\psi_\omega+o(\eps),\qquad
R_\omega=\frac{\partial\psi_\omega}{\partial \omega}+o(\eps).
\]
Therefore, since $\psi_\omega$ is radial,
\[
\int_{\R^N}y_j^2\varphi_\omega^2=\int_{\R^N}y_j^2\psi_\omega^2+o(\eps)=\frac{1}{N}\nldd{|y|\psi_\omega}+o(\eps),
\]
and so
\begin{multline}\label{eq:star}
\int_{\R^N} \nabla^2V(0)(y,y)\varphi_\omega T_\lambda=\\-\left(\frac{1}{p-1}-\frac{N+2}{4}\right)\frac{1}{N}\nldd{|y|\psi_\omega}\Delta V(0)+o(\eps).
\end{multline}
Similarly, we have 
\begin{multline}\label{eq:starstar}
\int_{\R^N}\nabla^2 Z(0)(y,y)R_\omega\varphi_\omega=\\-\left(\frac{1}{p-1}-\frac{N+2}{4}\right)
\frac{1}{N}\nldd{|y|\psi_\omega}\Delta Z(0)+o(\eps).
\end{multline}

Summarizing, using successively \eqref{eq:hessian1}, \eqref{eq:hessian2}, \eqref{eq:star}, \eqref{eq:starstar} and \eqref{eq:degencase} we have obtained
\begin{multline}\label{eq:summarize}
\int_{\R^N}R_\omega\varphi_\omega=-\frac{1}{2(\omega+V(0))}\nldd{\varphi_\omega}\\+\eps^2\frac{(\Delta Z(0)-\Delta V(0))}{NZ(0)}\left(\frac{1}{p-1}-\frac{N+2}{4}\right)\nldd{|y|\psi_\omega}+o(\eps^2).
\end{multline}

Now, we compute $\frac{\partial Q(\varphi_\omega)}{\partial\omega}$. 
First, recall that, coming back to the definition \eqref{eq:defQ} of $Q$, we have 
\[
\eps^{-N}\frac{\partial Q(\varphi_\omega)}{\partial\omega}=\nldd{\varphi_\omega}+2\omega\int_{\R^N}R_\omega\varphi_\omega+2\int_{\R^N}V(\eps y)R_\omega\varphi_\omega.
\]
As in \eqref{eq:Vexpansion}, \eqref{eq:cancelation}, and \eqref{eq:star} we can expand in $\eps$ and get 
\begin{multline*}
2\int_{\R^N}V(\eps y)R_\omega\varphi_\omega=2V(0)\int_{\R^N} R_\omega\varphi_\omega\\
-\eps^2\left(\frac{1}{p-1}-\frac{N+2}{4}\right)\frac{1}{N}\nldd{|y|\psi_\omega}\Delta V(0)+o(\eps^2).
\end{multline*}
Therefore, 
\begin{multline*}
\eps^{-N}\frac{\partial Q(\varphi_\omega)}{\partial\omega}=\nldd{\varphi_\omega}+2(\omega+V(0))\int_{\R^N}R_\omega\varphi_\omega\\
-\eps^2\left(\frac{1}{p-1}-\frac{N+2}{4}\right)\frac{1}{N}\nldd{|y|\psi_\omega}\Delta V(0)+o(\eps^2).
\end{multline*}
Using \eqref{eq:summarize}, we finally get
\begin{multline*}
\eps^{-N}\frac{\partial Q(\varphi_\omega)}{\partial\omega}=\eps^2\left(\frac{1}{p-1}-\frac{N+2}{4}\right)\frac{1}{N}\nldd{|y|\psi_\omega}\times\\
\times\left(\Delta Z(0)-\Delta V(0)\left(1+\frac{2(\omega+V(0))}{Z(0)}\right)\right)+o(\eps^2).
\end{multline*}

\subsection{The Spectral Condition}
We define the operator $L_0:=-\Delta+Z(0)-p\psi_\omega^{p-1}$. It is well known (see e.g. \cite{AmMa06}) that the 
spectrum of $L_0$ consists of one negative eigenvalue, a $N$-dimensional kernel (generated by $\frac{\partial \psi_\omega}{\partial y_j}$ for $j=1,\dots,N$) 
and the rest of the spectrum is bounded away from $0$. When $\eps$ is close to $0$, the spectrum of $L_\eps$ will be close to the spectrum of $L_0$. In particular, 
the $0$ eigenvalue, of multiplicity $N$, will transform into $N$ possibly different eigenvalues close to $0$ but shifted either to the positive or to the 
negative side of the real axis, depending on the sign of the eigenvalues of the Hessian of $Z$ at $0$. More precisely, the following proposition 
was proved in \cite{LiWe08} (see \cite{IaLe09} for a detailed justification).

\begin{proposition}\label{prop:spectrum}
The spectrum of $L_\eps$ consists of positive spectrum away from $0$ and a set of $N+1$ simple eigenvalues $\{\lambda_0,\lambda_1,\dots,\lambda_N\}$ such that 
\[
\lambda_0<\lambda_1\leq\cdots\leq\lambda_N.
\]
As $\eps\to0$, we have $\lambda_0<0$ and the following asymptotic expansion holds for the other eigenvalues:
\[
\lambda_j=c_j\eps^2+o(\eps^2),\qquad j=1,...,N,
\]
where $c_j=\frac{1}{2}\frac{\nldd{\psi_\omega}}{\nldd{\frac{\partial\psi_\omega}{\partial x_j}}}a_j$ and $\{a_1,\dots,a_N\}$ are the eigenvalues of the Hessian matrix $\nabla^2Z(0)$.
\end{proposition}
Therefore, (3) in Theorem \ref{mainthm} is a direct consequence of Proposition \ref{prop:spectrum}. In particular, the spectral condition for stability will be satisfied if and only if $0$ is a non-degenerate local minimum of $Z$.

\bibliographystyle{siam}
\bibliography{biblio}

\def\cprime{$'$}
\begin{thebibliography}{10}

\bibitem{AmBaCi97}
{\sc A.~Ambrosetti, M.~Badiale, and S.~Cingolani}, {\em Semiclassical states of
  nonlinear {S}chr\"odinger equations}, Arch. Rational Mech. Anal., 140 (1997),
  pp.~285--300.

\bibitem{AmMa06}
{\sc A.~Ambrosetti and A.~Malchiodi}, {\em Perturbation methods and semilinear
  elliptic problems on {${\bf R}^n$}}, vol.~240 of Progress in Mathematics,
  Birkh\"auser Verlag, Basel, 2006.

\bibitem{BS}
{\sc M.~Beals and W.~Strauss}, {\em {$L^p$} estimates for the wave equation
  with a potential}, Comm. Partial Differential Equations, 18 (1993),
  pp.~1365--1397.

\bibitem{BeCa81}
{\sc H.~Berestycki and T.~Cazenave}, {\em Instabilit\'e des \'etats
  stationnaires dans les \'equations de {S}chr\"odinger et de {K}lein-{G}ordon
  non lin\'eaires}, C. R. Acad. Sci. Paris S\'er. I Math., 293 (1981),
  pp.~489--492.

\bibitem{Ca03}
{\sc T.~Cazenave}, {\em Semilinear {S}chr\"odinger equations}, New York
  University -- Courant Institute, New York, 2003.

\bibitem{CaLi82}
{\sc T.~Cazenave and P.-L. Lions}, {\em Orbital stability of standing waves for
  some nonlinear {S}chr\"odinger equations}, Comm. Math. Phys., 85 (1982),
  pp.~549--561.

\bibitem{GrShSt87}
{\sc M.~Grillakis, J.~Shatah, and W.~A. Strauss}, {\em Stability theory of
  solitary waves in the presence of symmetry {I}}, J. Funct. Anal., 74 (1987),
  pp.~160--197.

\bibitem{GrShSt90}
\leavevmode\vrule height 2pt depth -1.6pt width 23pt, {\em Stability theory of
  solitary waves in the presence of symmetry {II}}, J. Funct. Anal., 94 (1990),
  pp.~308--348.

\bibitem{IaLe09}
{\sc I.~Ianni and S.~Le~Coz}, {\em Orbital stability of standing waves of a
  semiclassical nonlinear {S}chr\"odinger-{P}oisson equation}, Adv.
  Differential Equations, 14 (2009), pp.~717--748.

\bibitem{JeLe06}
{\sc L.~Jeanjean and S.~Le~Coz}, {\em An existence and stability result for
  standing waves of nonlinear {S}chr\"odinger equations}, Adv. Differential
  Equations, 11 (2006), pp.~813--840.

\bibitem{JeLe09}
\leavevmode\vrule height 2pt depth -1.6pt width 23pt, {\em Instability for
  standing waves of nonlinear {K}lein-{G}ordon equations via mountain-pass
  arguments}, Trans. Amer. Math. Soc., 361 (2009), pp.~5401--5416.

\bibitem{Le09}
{\sc S.~Le~Coz}, {\em Standing waves in nonlinear {S}chr\"odinger equations},
  in Analytical and numerical aspects of partial differential equations, Walter
  de Gruyter, Berlin, 2009, pp.~151--192.

\bibitem{LeFuFiKsSi08}
{\sc S.~Le~Coz, R.~Fukuizumi, G.~Fibich, B.~Ksherim, and Y.~Sivan}, {\em
  Instability of bound states of a nonlinear {S}chr\"odinger equation with a
  {D}irac potential}, Phys. D, 237 (2008), pp.~1103--1128.

\bibitem{LiWe08}
{\sc T.-C. Lin and J.~Wei}, {\em Orbital stability of bound states of
  semiclassical nonlinear {S}chr\"odinger equations with critical
  nonlinearity}, SIAM J. Math. Anal., 40 (2008), pp.~365--381.

\bibitem{LiOhTo07}
{\sc Y.~Liu, M.~Ohta, and G.~Todorova}, {\em Instabilit\'e forte d'ondes
  solitaires pour des \'equations de {K}lein-{G}ordon non lin\'eaires et des
  \'equations g\'en\'eralis\'ees de {B}oussinesq}, Ann. Inst. H. Poincar\'e
  Anal. Non Lin\'eaire, 24 (2007), pp.~539--548.

\bibitem{Oh89}
{\sc Y.-G. Oh}, {\em Stability of semiclassical bound states of nonlinear
  {S}chr\"odinger equations with potentials}, Comm. Math. Phys., 121 (1989),
  pp.~11--33.

\bibitem{OhTo05}
{\sc M.~Ohta and G.~Todorova}, {\em Strong instability of standing waves for
  nonlinear {K}lein-{G}ordon equations}, Discrete Contin. Dyn. Syst., 12
  (2005), pp.~315--322.

\bibitem{OhTo06}
\leavevmode\vrule height 2pt depth -1.6pt width 23pt, {\em Instability of
  standing waves for nonlinear {K}lein-{G}ordon equation and related system},
  in Nonlinear dispersive equations, vol.~26 of GAKUTO Internat. Ser. Math.
  Sci. Appl., Gakk\=otosho, Tokyo, 2006, pp.~189--200.

\bibitem{OhTo07}
\leavevmode\vrule height 2pt depth -1.6pt width 23pt, {\em Strong instability
  of standing waves for the nonlinear {K}lein-{G}ordon equation and the
  {K}lein-{G}ordon-{Z}akharov system}, SIAM J. Math. Anal., 38 (2007),
  pp.~1912--1931 (electronic).

\bibitem{Sh85}
{\sc J.~Shatah}, {\em Unstable ground state of nonlinear {K}lein-{G}ordon
  equations}, Trans. Amer. Math. Soc., 290 (1985), pp.~701--710.

\bibitem{ShSt85}
{\sc J.~Shatah and W.~A. Strauss}, {\em Instability of nonlinear bound states},
  Comm. Math. Phys., 100 (1985), pp.~173--190.

\bibitem{St08}
{\sc C.~A. Stuart}, {\em Lectures on the orbital stability of standing waves
  and application to the nonlinear {S}chr\"odinger equation}, Milan J. Math.,
  76 (2008), pp.~329--399.

\bibitem{Ta09}
{\sc T.~Tao}, {\em Why are solitons stable?}, Bull. Amer. Math. Soc., 46
  (2009), pp.~1--33.

\bibitem{We83}
{\sc M.~I. Weinstein}, {\em Nonlinear {S}chr\"odinger equations and sharp
  interpolation estimates}, Comm. Math. Phys., 87 (1982/83), pp.~567--576.

\bibitem{We85}
{\sc M.~I. Weinstein}, {\em Modulational stability of ground states of
  nonlinear {S}chr\"odinger equations}, SIAM J. Math. Anal., 16 (1985),
  pp.~472--491.

\end{thebibliography}

%
\end{document}